\documentclass[12pt]{article}
\usepackage[T1]{fontenc}
\usepackage[utf8]{inputenc}
\usepackage{lmodern}
\usepackage[a4paper,margin=27mm]{geometry}
\usepackage{amsmath,amssymb,amsthm,mathtools}
\usepackage{booktabs,array,enumitem,microtype}
\usepackage{xcolor}
\usepackage[colorlinks=true,linkcolor=blue!45!black,citecolor=blue!45!black,urlcolor=blue!45!black]{hyperref}
\usepackage{fancyhdr}
\usepackage{authblk}
\numberwithin{equation}{section}
\newtheorem{theorem}{Theorem}[section]
\newtheorem{lemma}[theorem]{Lemma}
\newtheorem{proposition}[theorem]{Proposition}
\newtheorem{corollary}[theorem]{Corollary}
\theoremstyle{remark}
\newcommand{\R}{\mathbb R}
\newcommand{\Z}{\mathbb Z}
\newcommand{\SO}{\mathrm{SO}}
\newcommand{\OO}{\mathrm O}
\newcommand{\SU}{\mathrm{SU}}
\newcommand{\Sp}{\mathrm{Sp}}
\newcommand{\Spin}{\mathrm{Spin}}
\newcommand{\Isom}{\operatorname{Isom}}

\newcommand{\diag}{\operatorname{diag}}
\newcommand{\tr}{\operatorname{tr}}
\newcommand{\rank}{\operatorname{rank}}
\newcommand{\Lie}{\operatorname{Lie}}
\newcommand{\Span}{\operatorname{span}}
\newcommand{\g}{\mathfrak g}
\newcommand{\h}{\mathfrak h}
\newcommand{\m}{\mathfrak m}
\newcommand{\dd}{\mathbf d}
\newcommand{\Md}{M_{\dd}}
\newcommand{\cohom}{\operatorname{cohom}}
\newcommand{\doi}[1]{\href{https://doi.org/#1}{\nolinkurl{#1}}}
\newcommand{\arxiv}[1]{\href{https://arxiv.org/abs/#1}{arXiv:\nolinkurl{#1}}}
\hypersetup{pdftitle={Normal Homogeneous Counterexamples to Wolf's Homogeneity Conjecture},pdfsubject={Revised manuscript on normal homogeneous Clifford-Wolf quotients},pdfkeywords={Clifford-Wolf subgroup, homogeneity conjecture, normal homogeneous metric, real flag manifold}}
\title{Normal Homogeneous Counterexamples to\\ Wolf's Homogeneity Conjecture}

\author[1]{Zaili Yan}
\author[2]{Hui Zhang}
\affil[1]{School of Mathematics and Statistics, Ningbo University, Ningbo, Zhejiang Province, 315211, People's Republic of China, yanzaili@nbu.edu.cn}
\affil[2]{School of Mathematics, Southeast University, Nanjing 211189, People's Republic of China, huizhang@mail.nankai.edu.cn}
\date{}
\begin{document}
\maketitle

\begin{abstract}
We construct free Clifford--Wolf subgroups of order two whose quotients of compact, connected, simply connected normal homogeneous manifolds are not homogeneous. The universal covers are the oriented real flag manifolds
\(M_{\dd}=\SO(N)/\prod_{i=1}^{q}\SO(d_i)\), equipped with the metric induced by the negative Killing form,
 where $N\ge4$ and $q\ge3$ is odd.  We determine its full connected centralizer, the effective connected isometry group of the quotient, and the corresponding orbit space. These calculations prove nonhomogeneity and give explicit formulas for the cohomogeneity, fundamental group, and Euler characteristic.
 These examples are found by ChatGPT 6 Astra.
\end{abstract}

\section{Introduction and the main theorem}

An isometry $f$ of a connected Riemannian manifold $(M,g)$ is a \emph{Clifford--Wolf transformation}, or a \emph{CW transformation}, if its displacement function
\[
\delta_f(x)=d_g(x,f(x))
\]
is constant on $M$. A subgroup $\Gamma\le\Isom(M,g)$ is a \emph{CW subgroup} if every one of its elements is a CW transformation. These definitions concern the distance on the entire manifold; they do not merely impose a condition on the differential of an isometry at one point.

Wolf's homogeneity problem asks about the passage from a homogeneous universal Riemannian covering to a globally homogeneous quotient. In its usual formulation, $M$ is connected, simply connected and homogeneous, and $\Gamma$ acts freely and properly discontinuously by isometries. The proposed equivalence is
\begin{equation}\label{eq:conjecture}
\Gamma\backslash M\text{ is homogeneous}
\quad\Longleftrightarrow\quad
\Gamma\text{ is a CW subgroup of }\Isom(M,g).
\end{equation}
The forward implication is classical \cite[Theorems~1 and~2]{WolfCover}. Indeed, the connected transitive isometry group of the quotient lifts to a connected transitive group on $M$ that centralizes the deck group. If $a$ belongs to this lifted group and $\gamma\in\Gamma$, then
\[
 d_g(a(x),\gamma a(x))=d_g(a(x),a\gamma (x))=d_g(x,\gamma (x)).
\]
Transitivity makes the displacement independent of $x$.

The converse in \eqref{eq:conjecture}, traditionally called \emph{Wolf's Homogeneity Conjecture} \cite{WolfSurvey},
 has been studied for several decades.
It holds when $M$ is Riemannian symmetric \cite{Freudenthal,Ozols,WolfVincent,WolfSym},
 has non-positive sectional curvature \cite{Druetta,WolfNegative},
 admits an invariant metric of strictly positive sectional curvature \cite{WolfPositive,WolfSpheres},
  a noncompact semisimple transitive group \cite{DottiMiatelloWolf},  or an exponential solvable transitive group \cite{WolfBounded,WolfCorrection}.
In \cite{DengWolf}, Deng and Wolf considered an analog of the Homogeneity Conjecture in Finsler
geometry, and proved it for symmetric Finsler spaces.
See Wolf's survey \cite{WolfSurvey} for more information.
In general, the Homogeneity Conjecture is false. Xu and Deng constructed left-invariant metrics on $\Sp(2)$ and cyclic CW subgroups of odd order whose quotients are not homogeneous \cite[Theorem~1.2]{XuDengCounter}.

We mention that CW subgroups should also be distinguished from \emph{CW homogeneity}, which requires CW transformations carrying any prescribed point to any other. The latter condition is considerably stronger and has its own classification \cite{BerestovskiiNikonorov}. Likewise, a constant-length Killing vector field supplies infinitesimal information; a finite CW transformation need not lie in an all-time CW flow. The classification of certain constant-length Killing vector fields on normal homogeneous spaces in \cite{XuWolf} therefore does not replace the distance calculation needed here.

We now introduce the main results in this paper. Fix an ordered orthogonal decomposition
\[
 \R^N=V_1\oplus\cdots\oplus V_q,
 \qquad d_i=\dim V_i\ge1,
 \qquad N=\sum_{i=1}^q d_i,
 \qquad q\text{ odd}.
\]
Write $\dd=(d_1,\ldots,d_q)$ and set
\[
 H=\prod_{i=1}^q\SO(V_i)\subset\SO(N),
 \qquad \Md=\SO(N)/H.
\]
We interpret $\SO(1)$ as the trivial group. A point of $\Md$ is an ordered orthogonal decomposition into oriented subspaces of dimensions $d_i$, with compatible ambient orientation. A singleton block records an oriented line, equivalently a unit vector. We call $\Md$ an \emph{oriented real flag manifold}; this terminology does not assert that it is a complex generalized flag manifold.

Let $\g=\mathfrak{so}(N)$ and $\h=\bigoplus_i\mathfrak{so}(V_i)$ be the Lie algebras of $\SO(N)$ and $H$ respectively. The positive-definite form
\[
 Q(X,Y)=-K_{\mathfrak{so}(N)}(X,Y)=-(N-2)\tr(XY)
 \qquad(X,Y\in\g)
\]
defines a bi-invariant metric on $\SO(N)$. Put $\m=\h^{\perp_Q}$ the orthogonal complement of $\h$ in $\g$ with respect to $Q$. The quotient metric $g$ is defined by identifying $T_{eH}\Md$ with $(\m,Q|_{\m})$ and translating by $\SO(N)$. Equivalently, $\SO(N)\to\Md$ is a Riemannian submersion. This is the normal metric used throughout; all displacement formulas refer to this normalization.

Choose a unit vector $u\in\R^N$ and let
\[
 r_u=I_N-2uu^T.
\]
Thus $r_uu=-u$, $r_u|_{u^\perp}=I$, and $\det r_u=-1$. Choose
\begin{equation}\label{eq:w}
 w=\diag(w_1,\ldots,w_q),\qquad
 w_i\in\OO(V_i),\qquad \det w_i=-1\quad\text{for every }i.
\end{equation}
In particular, $\det w=(-1)^q=-1$ and $w\in N_{\OO(N)}(H)$. Define
\begin{equation}\label{eq:f}
 f_{u,w}(xH)=r_uxwH.
\end{equation}
The product $r_uxw$ belongs to $\SO(N)$. The individual block conditions in \eqref{eq:w} are essential. Any two choices satisfying \eqref{eq:w} differ by right multiplication by an element of $H$, so the induced map is independent of that choice. We retain the notation $f_{u,w}$ because it is convenient for the matrix calculations.

\begin{theorem}\label{thm:main}
Assume that $N\ge4$, that $q\ge3$ is odd, and that $d_i\ge2$ for at least one $i$. Equip $\Md=\SO(N)/H$ with the normal metric $g$ defined above, and choose $w$ as in \eqref{eq:w}.
\begin{enumerate}[label=(\roman*),leftmargin=2.1em]
\item The manifold $\Md$ is compact, connected and simply connected. The map $f_{u,w}$ is a CW involution of displacement
\begin{equation}\label{eq:delta}
 \delta_{f_{u,w}}=\pi\sqrt{(N-2)(q-1)}.
\end{equation}
Consequently $\Gamma=\langle f_{u,w}\rangle\cong\Z_2$ acts freely and is a CW subgroup.
\item The Riemannian quotient $\Gamma\backslash\Md$ is not homogeneous. Let
\[
 V_0=\bigoplus_{d_i=1}V_i,\qquad t=\dim V_0=\#\{i:d_i=1\}.
\]
The cohomogeneity of its full connected isometry group is
\[
 \cohom\bigl(\Isom(\Gamma\backslash\Md)^0,\Gamma\backslash\Md\bigr)
 =\begin{cases}q-1,&t=0,\\q-t,&t>0.\end{cases}
\]
Here cohomogeneity means the dimension of the orbit space, or equivalently the codimension of a principal orbit.
\end{enumerate}
\end{theorem}

This paper is organized as follows. Section~2 proves Theorem \ref{thm:main} (i) by a sharp lower bound for the multiplicity of the eigenvalue $-1$.
 Section~3 determines the connected isometry group and centralizer, proves Theorem \ref{thm:main} (ii), and computes the fundamental group and Euler characteristic.
 Section~4 specializes the construction to Stiefel and full flag manifolds.

\section{Normal quotient distance and proof of Theorem 1.1 (i)}\label{sec:distance}
We use $d_Q$ for the distance on $\SO(N)$ associated with the bi-invariant metric $Q$ and $d_g$ for the quotient distance on $\Md$.

\begin{lemma}[\cite{DolcettiPertici}]\label{lem:angles}
Let $A\in\SO(N)$ have principal rotation angles
$\theta_1,\ldots,\theta_{\lfloor N/2\rfloor}\in[0,\pi]$, including zero angles when necessary. Then
\begin{equation}\label{eq:angles}
 d_Q(I_N,A)^2=2(N-2)\sum_{j=1}^{\lfloor N/2\rfloor}\theta_j^2.
\end{equation}
\end{lemma}

\begin{proposition}\label{prop:minus}
Let $A\in\SO(N)$ and let $a$ be a nonnegative integer. If
$\dim\ker(A+I_N)\ge2a$, then
\[
 d_Q(I_N,A)^2\ge2a(N-2)\pi^2.
\]
Equality holds if and only if $A$ is an involution with exactly $2a$ negative eigenvalues.
\end{proposition}
\begin{proof}
The multiplicity of $-1$ in an element of $\SO(N)$ is even. In its real orthogonal normal form, each pair of negative eigenvalues gives a rotation block of angle $\pi$. Thus at least $a$ of the angles in \eqref{eq:angles} equal $\pi$, and every remaining summand is nonnegative. This proves the inequality. Equality requires exactly $a$ angles equal to $\pi$ and every other angle equal to zero. The eigenvalues are then precisely $2a$ copies of $-1$ and $N-2a$ copies of $1$. Conversely, that spectrum for an orthogonal matrix gives the asserted equality. The case $a=0$ says that equality occurs exactly at the identity.
\end{proof}

The Riemannian-submersion distance formula is
\begin{equation}\label{eq:quotientdistance}
 d_g(xH,yH)=\min_{h\in H}d_Q(x,yh)
           =\min_{h\in H}d_Q(I_N,x^{-1}yh).
\end{equation}
This is the standard normal quotient formula; see \cite{ONeill} and \cite[Chapter~7]{Besse}.
Compactness of $H$ ensures that the minimum in \eqref{eq:quotientdistance} is attained.

\begin{proposition}\label{prop:reduction}
For $x\in\SO(N)$ put $v=x^{-1}u$. Then
\begin{equation}\label{eq:reduction}
 d_g(xH,f_{u,w}(xH))
 =\min_{h\in H}d_Q(I_N,r_vwh).
\end{equation}
\end{proposition}
\begin{proof}
Orthogonality of $x$ gives $\|v\|=1$ and
\[
 x^{-1}r_ux
 =x^{-1}(I_N-2uu^T)x
 =I_N-2(x^{-1}u)(x^{-1}u)^T=r_v.
\]
Substitution into \eqref{eq:quotientdistance} yields
\[
\begin{aligned}
 d_g(xH,r_uxwH)
 &=\min_{h\in H}d_Q(I_N,x^{-1}r_uxwh)\\
 &=\min_{h\in H}d_Q(I_N,r_vwh),
\end{aligned}
\]
as required.
\end{proof}

\textbf{Proof of Theorem~\ref{thm:main} (i).}
We separate the topological, isometric and distance arguments.

\emph{Topology.} The group $H$ is closed and connected, and $\SO(N)$ is compact and connected. Hence $\Md$ is compact and connected. The homotopy sequence of the principal bundle $H\to\SO(N)\to\Md$ contains
\[
 \pi_1(H)\longrightarrow\pi_1(\SO(N))
 \longrightarrow\pi_1(\Md)\longrightarrow\pi_0(H)=0.
\]
For $N\ge3$, $\pi_1(\SO(N))\cong\Z_2$, generated by one full rotation in a coordinate two-plane. Since at least one $d_i$ is at least two, this loop lies in $\SO(V_i)\subset H$. The first map is therefore surjective, including when $N=4$. Exactness gives $\pi_1(\Md)=0$. $\Md$ is simply connected.

\emph{The isometry property.}  On $\SO(N)$ consider $F(x)=r_uxw$. For a tangent vector $xX$, $X\in\mathfrak{so}(N)$,
\[
 dF_x(xX)=r_uxXw=F(x)(w^{-1}Xw).
\]
Conjugation by $w\in\OO(N)$ preserves the trace form:
\[
 Q(w^{-1}Xw,w^{-1}Yw)
 =-(N-2)\tr(w^{-1}XYw)=Q(X,Y).
\]
Thus $F$ is an isometry of $\SO(N)$. It maps right $H$-fibers to right $H$-fibers, so it induces an isometry of the quotient. Since $w_i^2\in\SO(V_i)$ on every block, we have
\[
 w^2\in H,\qquad
 f_{u,w}^{2}(xH)=r_u(r_uxw)wH=xw^2H=xH.
\]
Thus the  induced isometry $f_{u,w}$ is therefore an involution.

\emph{The lower bound.} Fix $x$ and put $v=x^{-1}u$. For $h=\diag(h_1,\ldots,h_q)\in H$, let
\[
 S=wh=\diag(S_1,\ldots,S_q),\qquad S_i=w_ih_i.
\]
Every $S_i$ is orthogonal with determinant $-1$. Its nonreal eigenvalues occur in conjugate pairs with product one, and its real eigenvalues are $1$ or $-1$. Hence the number of negative eigenvalues is odd, in particular positive. The $q$ blocks give
\[
 E:=\ker(S+I_N),\qquad \dim E\ge q.
\]
Intersecting with one hyperplane can reduce dimension by at most one:
\[
 \dim(E\cap v^\perp)\ge\dim E-1\ge q-1.
\]
If $z\in E\cap v^\perp$, then $Sz=-z$ and $r_vz=z$, so
\[
 r_vSz=-z,
 \qquad E\cap v^\perp\subset\ker(r_vS+I_N).
\]
The product $r_vS$ belongs to $\SO(N)$. Because $q$ is odd, $(q-1)/2$ is an integer, and Proposition~\ref{prop:minus} gives
\begin{equation}\label{eq:lower}
 d_Q(I_N,r_vS)^2\ge(N-2)(q-1)\pi^2
 \qquad(S\in wH).
\end{equation}
Taking the minimum in \eqref{eq:reduction} proves the same lower bound for the squared quotient displacement.

\emph{An explicit representative attaining the bound.} Decompose
\[
 v=v_1+\cdots+v_q,\qquad v_i\in V_i.
\]
If $v_i\ne0$, set $n_i=v_i/\|v_i\|$; if $v_i=0$, choose any unit vector $n_i\in V_i$. Define
\[
 S_i=I_{V_i}-2n_in_i^T,\qquad
 S=\diag(S_1,\ldots,S_q).
\]
Then $\det S_i=-1$ and $h_i=w_i^{-1}S_i\in\SO(V_i)$, so $S=wh$ for $h=\diag(h_1,\ldots,h_q)\in H$. The vectors $n_i$ are orthonormal because they belong to mutually orthogonal blocks. Let
\[
 W=\Span\{n_1,\ldots,n_q\},\qquad \dim W=q.
\]
We have $S|_W=-I_W$, $S|_{W^\perp}=I$, and $v\in W$. With respect to the orthogonal decomposition
\[
 \R^N=\R v\oplus(W\cap v^\perp)\oplus W^\perp,
\]
the two reflections have the forms
\[
 r_v=\diag(-1,I_{q-1},I_{N-q}),\qquad
 S=\diag(-1,-I_{q-1},I_{N-q}),
\]
and consequently
\[
 r_vS=\diag(1,-I_{q-1},I_{N-q}).
\]
This is an involution with exactly $q-1$ negative eigenvalues. Proposition~\ref{prop:minus} gives
\[
 d_Q(I_N,r_vS)^2=(N-2)(q-1)\pi^2.
\]
Together with \eqref{eq:lower}, this proves \eqref{eq:delta} at every point $xH$.

Finally, $N\ge4$ and $q\ge3$ make \eqref{eq:delta} positive. Therefore $f_{u,w}\ne1$ and it has no fixed point. Since its square is the identity, $\Gamma=\{1,f_{u,w}\}\cong\Z_2$ is free. Both of its elements have constant displacement, so it is a CW subgroup. A finite free isometric action is properly discontinuous, and the quotient is a smooth Riemannian manifold.
\qed

\section{The full connected centralizer and nonhomogeneity}\label{sec:centralizer}
For a connected compact Riemannian manifold $M$, write $I=\Isom(M)$ and $I^0$ for its identity component. The Myers--Steenrod theorem makes $I$ a compact Lie group \cite{MyersSteenrod}. For a subgroup $A\le I$, denote its normalizer and centralizer by $N_I(A)$ and $C_I(A)$.
\begin{proposition}[\cite{WolfSym}]\label{prop:cover}
Let $M$ be connected, compact and simply connected, and let $\Gamma\le I$ be finite and free. Then
\begin{equation}\label{eq:covergroup}
 \Isom(\Gamma\backslash M)\cong N_I(\Gamma)/\Gamma,
 \qquad
 \Isom(\Gamma\backslash M)^0
 \cong C_I(\Gamma)^0/(C_I(\Gamma)^0\cap\Gamma).
\end{equation}
In particular, $\Gamma\backslash M$ is homogeneous if and only if $C_I(\Gamma)^0$ is transitive on $M$.
\end{proposition}

We now return to $\Md$. Rotations of $V_0$ commute with $H$, normalize it, and act isometrically from the right. We use the convention
\[
 R_c(xH)=xc^{-1}H,\qquad c\in\SO(V_0),
\]
with $c$ extended by the identity on all nonsingleton blocks. We set $\SO(0)=\SO(1)=\{1\}$ when these groups occur in formulas. Define
\[
 \Phi:\mathrm{SO}(N)\times\mathrm{SO}(V_0)\longrightarrow I^0,
 \qquad \Phi(a,c)=L_aR_c.
\]
\begin{lemma}\label{lemkernel}
 The kernel of $\Phi$ is given by
 \begin{equation*}
 \ker\Phi=
 \begin{cases}
 \{1,(-I_{N},-I_{t})\},
     &t\, \text{even}, d_i\, \text{even if }\, d_i\geq2,\\
 \{1\},&\text{otherwise}.
 \end{cases}
\end{equation*}
\end{lemma}
\begin{proof}
The left action has discrete kernel: for $N\ge5$ this follows from simplicity and $\mathfrak h\ne\mathfrak{so}(N)$, while for $N=4$ the one-dimensional algebra $\mathfrak h$ contains no nonzero ideal of $\mathfrak{so}(4)$. If $\Phi(a,c)=1$, then $L_a=R_c^{-1}$ commutes with every left translation. Differentiating this equality shows that $\operatorname{Ad}(a)$ is the identity on $\mathfrak{so}(N)$, hence $a=\varepsilon I_N$ for some $\varepsilon\in\{1,-1\}$. Evaluation at $eH$ gives $ac^{-1}\in H$. On the singleton subspace this forces $c=\varepsilon I_t$, while on every nonsingleton block it requires $\varepsilon^{d_i}=1$. Conversely, these conditions suffice. Thus
\[
 \ker\Phi
 =\left\{
 (\varepsilon I_N,\varepsilon I_t):
 \begin{array}{l}
 \varepsilon\in\{1,-1\},\quad\varepsilon^N=\varepsilon^t=1,\\
 \varepsilon^{d_i}=1\text{ for every }i\text{ with }d_i\ge2
 \end{array}
 \right\}.
\]
 In particular, the only possible nonidentity kernel element is $(-I_N,-I_t)$; it occurs precisely when every nonsingleton block dimension and $t$ are even.
\end{proof}

\begin{proposition}\label{lem:isometry}
Under the hypotheses of Theorem~\ref{thm:main}, $(\Md,g)$ is locally irreducible and has nonconstant sectional curvature. Its full connected isometry group is the effective image
\begin{equation}\label{eq:connected}
 \Isom(\Md,g)^0=L(\SO(N))R(\SO(V_0))/\ker \Phi.
\end{equation}
\end{proposition}
\begin{proof}
We first verify local irreducibility, including the low-dimensional case $N=4$.
Choose an orthonormal basis adapted to the ordered block decomposition and put
\[
 E_{ab}=e_ae_b^{T}-e_be_a^{T}\qquad(a\ne b).
\]
The space $\mathfrak m=\mathfrak h^{\perp_Q}$ is spanned by the matrices $E_{ab}$ whose indices belong to different blocks. For indices $a,b$ in the same block, choose an index $c$ in a different block. Then
\[
 E_{ac},E_{bc}\in\mathfrak m,
 \qquad [E_{ac},E_{bc}]=-E_{ab}.
\]
Thus $\mathfrak m$, together with its brackets, contains every elementary skew-symmetric matrix, and consequently
\[
 \langle\mathfrak m\rangle_{\mathrm{Lie}}=\mathfrak{so}(N).
\]
Suppose that the de Rham decomposition of $M_{\mathbf d}$ is nontrivial. The connected group of left translations preserves the factor distributions. Grouping the factors into two nonempty collections gives an orthogonal splitting
\[
 \mathfrak m=\mathfrak m_1\oplus\mathfrak m_2,
 \qquad \mathfrak m_1\ne0,\quad\mathfrak m_2\ne0.
\]
Every plane spanned by $X\in\mathfrak m_1$ and $Y\in\mathfrak m_2$ has zero sectional curvature. The normal homogeneous curvature formula \cite[Chapter~7]{Besse}, obtained from the submersion formulas in \cite{ONeill}, therefore implies
\[
 0=Q(R(X,Y)Y,X)
   =\|[X,Y]_{\mathfrak h}\|_Q^2
      +\frac14\|[X,Y]_{\mathfrak m}\|_Q^2,
 \qquad [\mathfrak m_1,\mathfrak m_2]=0.
\]
Let $\mathfrak g_i$ be the Lie subalgebra generated by $\mathfrak m_i$, $i=1,2$. Repeated application of the Jacobi identity gives $[\mathfrak g_1,\mathfrak g_2]=0$. Since $\mathfrak m$ generates $\mathfrak{so}(N)$, it follows that
\[
 \mathfrak{so}(N)=\mathfrak g_1+\mathfrak g_2.
\]
Both summands are ideals, and their intersection is contained in the centre of $\mathfrak{so}(N)$. As this centre is zero, the sum is direct. For $N\ge5$, this contradicts the simplicity of $\mathfrak{so}(N)$.

If $N=4$, the assumptions force the block dimensions to be $(2,1,1)$, up to permutation. The only two nonzero proper ideals of $\mathfrak{so}(4)$ are its self-dual and anti-self-dual factors $\mathfrak{so}(3)_+$ and $\mathfrak{so}(3)_-$. They are orthogonal for $Q$. Hence the preceding decomposition would give
\[
 \mathfrak h
 =\bigl(\mathfrak g_1\cap\mathfrak m_1^{\perp_Q}\bigr)
  \oplus
  \bigl(\mathfrak g_2\cap\mathfrak m_2^{\perp_Q}\bigr).
\]
Since $\dim\mathfrak h=1$, the algebra $\mathfrak h$ would have to lie in one of these two ideals. After relabelling the basis, however, $\mathfrak h=\mathbb R E_{12}$, whereas
\[
 E_{12}
   =\frac12(E_{12}+E_{34})
     +\frac12(E_{12}-E_{34})
\]
has nonzero components in both ideals. This is a contradiction. Local irreducibility follows in every case.

We next exclude a round metric by an explicit curvature computation. Choose distinct indices $a,b$ in a block of dimension at least two, and choose indices $c,d$ in two other blocks. The matrices $E_{ac},E_{bd}\in\mathfrak m$ commute, so
\[
 K(\operatorname{span}\{E_{ac},E_{bd}\})=0.
\]
On the other hand, $[E_{ac},E_{bc}]=-E_{ab}\in\mathfrak h$, and
\[
 Q(E_{ij},E_{ij})=2(N-2),\qquad
 K(\operatorname{span}\{E_{ac},E_{bc}\})
 =\frac{2(N-2)}{\bigl(2(N-2)\bigr)^2}
 =\frac{1}{2(N-2)}>0.
\]
The sectional curvature is therefore not constant. In particular, the manifold is not isometric to a round sphere or a round real projective space.

Reggiani's connected-isometry theorem \cite[Corollary~1.3]{Reggiani} now applies: the identity component of the isometry group is the almost direct product of the semisimple transitive group and the connected group generated by its invariant Killing vector fields. To identify the latter group, use the block-module decomposition
\[
 \mathfrak m=\bigoplus_{i<j}V_i\wedge V_j.
\]
If $d_i\ge2$, the standard representation of $\mathrm{SO}(V_i)$ has no nonzero fixed vector. Since the factors of $H$ act independently, $V_i\wedge V_j$ has no $H$-fixed vector whenever at least one of $d_i,d_j$ exceeds one. The only fixed summands correspond to pairs of singleton blocks. Consequently
\[
 \mathfrak m^H
 =\bigwedge\nolimits^2 V_0
 =\mathfrak{so}(V_0).
\]
Evaluation at the base point identifies invariant vector fields with $\mathfrak m^H$. The right rotations $R_c$, $c\in\mathrm{SO}(V_0)$, commute with all left translations and realize exactly these fields. Their action is effective, because $\mathrm{SO}(V_0)\cap H=\{I_N\}$. Thus
\[
 I^0=L(\mathrm{SO}(N))R(\mathrm{SO}(V_0))/\ker\Phi,
\]
as an effective image, as asserted.
\end{proof}

\begin{proposition}\label{prop:centralizer}
Let $\Gamma=\langle f_{u,w}\rangle$ and $C=C_{\Isom(\Md,g)}(f_{u,w})^0$. Then
\[
 C\cong\SO(N-1)\times\SO(t),
\]
where $\SO(N-1)$ fixes $u$ and acts on $u^\perp$. Moreover,
\begin{equation}\label{eq:quotientgroup}
 \Isom(\Gamma\backslash\Md)^0\cong
 \begin{cases}
 \bigl(\SO(N-1)/\{\pm I_{N-1}\}\bigr)\times\SO(t),
     &\text{all }d_i\text{ odd},\\
 \SO(N-1)\times\SO(t),&\text{otherwise}.
 \end{cases}
\end{equation}
\end{proposition}
\begin{proof}
Write $f=f_{u,w}$ and let
\[
 K_u=\{a\in\mathrm{SO}(N):au=u\}\cong\mathrm{SO}(N-1).
\]
Because every singleton block of $w$ equals $-1$, we have $w|_{V_0}=-I_{V_0}$. Hence $w$ commutes with every $c\in\mathrm{SO}(V_0)$. Using the inverse lift $x\mapsto r_uxw^{-1}$, direct substitution yields
\[
 \begin{aligned}
 f\Phi(a,c)f^{-1}(xH)
  &=r_u a r_u xw^{-1}c^{-1}wH\\
  &=r_u a r_u xc^{-1}H
   =\Phi(r_u a r_u,c)(xH).
 \end{aligned}
\]
By Lemma \ref{lemkernel}, the kernel of $\Phi$ is discrete. Therefore the Lie algebra of the centralizer of $f$ is the fixed subalgebra of
\[
 (A,B)\longmapsto(r_u A r_u,B)
 \quad\text{on}\quad
 \mathfrak{so}(N)\oplus\mathfrak{so}(t).
\]
Relative to $\mathbb R u\oplus u^{\perp}$, write
\[
 A=\begin{pmatrix}0&-b^T\\b&D\end{pmatrix},
 \qquad
 r_u A r_u=\begin{pmatrix}0&b^T\\-b&D\end{pmatrix}.
\]
Thus $r_u A r_u=A$ if and only if $b=0$, and consequently
\[
 \operatorname{Lie}(C_I(f))
   =\mathfrak{so}(u^{\perp})\oplus\mathfrak{so}(t).
\]
The connected subgroup $\Phi(K_u\times\mathrm{SO}(t))$ has exactly this Lie algebra. It follows that
\[
 C_I(f)^0=\Phi(K_u\times\mathrm{SO}(t)).
\]
Moreover, the restriction of $\Phi$ to $K_u\times\mathrm{SO}(t)$ is injective: its only possible nonidentity kernel element has first coordinate $-I_N$, which does not fix $u$. Therefore
\[
 C_I(f)^0\cong\mathrm{SO}(N-1)\times\mathrm{SO}(t).
\]

It remains to compute $C_I(f)^0\cap\Gamma$. Since $\Gamma=\{1,f\}$, the only question is whether $f$ belongs to this connected centralizer. Suppose that $f=\Phi(a,c)$ with $a\in K_u$. Comparing conjugation on the left-translation Lie algebra gives
\[
 \operatorname{Ad}(a)|_{\mathfrak{so}(N)}
   =\operatorname{Ad}(r_u)|_{\mathfrak{so}(N)}.
\]
The commutant of the standard $\mathfrak{so}(N)$ representation in $O(N)$ is $\{\pm I_N\}$, so $a=\pm r_u$. The condition $au=u$ forces
\[
 a=-r_u,
 \qquad \det(-r_u)=(-1)^{N+1}.
\]
In particular, $N$ must be odd. Under this condition,
\[
 f(xH)=(-r_u)x(-w)H.
\]
Equality with $\Phi(-r_u,c)$ is equivalent to $(-w)c\in H$. On each nonsingleton block this says
\[
 \det(-w_i)=(-1)^{d_i}\det(w_i)=(-1)^{d_i+1}=1,
\]
so every such $d_i$ must be odd. On the singleton subspace, $-w=I_{V_0}$ and $H$ acts trivially, so necessarily $c=I_t$. Conversely, if every $d_i$ is odd, then $N$ is odd because $q$ is odd, $-w\in H$, and
\[
 f=L_{-r_u}\in C_I(f)^0.
\]
On $u^{\perp}$, the matrix $-r_u$ is $-I_{N-1}$. Accordingly,
\[
 C_I(f)^0\cap\Gamma
 =\begin{cases}
 \langle(-I_{N-1},I_t)\rangle,&\text{if every }d_i\text{ is odd},\\
 \{1\},&\text{otherwise}.
 \end{cases}
\]
In the first case $N-1$ is even, so $-I_{N-1}\in\mathrm{SO}(N-1)$. Applying the covering-space formula of Proposition~\ref{prop:cover} gives precisely the asserted two alternatives for $\operatorname{Isom}(\Gamma\backslash M_{\mathbf d})^0$.
\end{proof}

Let
\[
 \Delta^s=\{(p_0,\ldots,p_s)\in\R^{s+1}:
                 p_j\ge0,\ \textstyle\sum_{j=0}^{s}p_j=1\}
\]
be the closed $s$-simplex.
\begin{theorem}\label{thm:orbits}
The orbit space of the full connected isometry group of $\Gamma\backslash\Md$ is homeomorphic to
\begin{equation}\label{eq:orbits}
 \begin{cases}
 \Delta^{q-1},&t=0,\\
 \bigl\{(p_1,\ldots,p_{q-1},z)\in\R^{q}:p_i\ge0,
                     \ \sum_{i=1}^{q-1} p_i+z^2=1\bigr\},&t=1,\\
 \Delta^{q-t},&t\ge2.
 \end{cases}
\end{equation}
Its cohomogeneity is $q-1$ when $t=0$ and $q-t$ when $t>0$.
\end{theorem}
\begin{proof}
Put $K_u=\{a\in\SO(N):au=u\}$ and $C=C_{\Isom(\Md,g)}(f_{u,w})^0$.
 By Proposition~\ref{prop:centralizer}, $C$ is the image of $K_u\times \SO(V_0)$. The map
\[
 \SO(N)\longrightarrow S^{N-1},\qquad x\longmapsto v=x^{-1}u,
\]
identifies $K_u\backslash\SO(N)$ with the unit sphere.  Hence
\[
 C\backslash\Md\cong(H\times \SO(V_0))\backslash S^{N-1},
\]
where the action on the right is the orthogonal block action.

Let $J=\{i:d_i\ge2\}$ and $r=|J|=q-t$. Write
\[
 v=\sum_{i\in J}v_i+v_0,
 \qquad v_i\in V_i,\quad v_0\in V_0,
 \qquad \sum_{i\in J}\|v_i\|^2+\|v_0\|^2=1.
\]
For every $d_i\ge2$, the group $\SO(V_i)$ is transitive on each sphere of positive radius in $V_i$, and the origin is one orbit. Its rotations are independent of those in the other blocks. It follows that the following lists of quantities are complete invariants.

If $t=0$, take $p_i=\|v_i\|^2$ for $1\le i\le q$. They are nonnegative and sum to one. Conversely, every such list is realized by choosing $v_i=\sqrt{p_i}\,n_i$ for fixed unit vectors $n_i\in V_i$. Two vectors with the same list are related by an element of $H$. The orbit space is therefore $\Delta^{q-1}$.

If $t=1$, write $v_0=ze_0$, $e_0$ is an oriented unit vector. The group $\SO(1)$ is trivial, so the signed coordinate $z$, not only $z^2$, must be retained. The complete invariants are
\[
 (\|v_1\|^2,\ldots,\|v_{q-1}\|^2,z),
 \qquad \sum_{i=1}^{q-1}\|v_i\|^2+z^2=1.
\]
This gives the middle space in \eqref{eq:orbits}.

If $t\ge2$, the group $\SO(V_0)$ is transitive on spheres in $V_0$. The complete invariants are
\[
 p_0=\|v_0\|^2,\qquad p_i=\|v_i\|^2\ (i\in J),
 \qquad p_0+\sum_{i\in J}p_i=1,
\]
which parametrize $\Delta^r=\Delta^{q-t}$. In every case, the invariant map is continuous, surjective and separates orbits. Since its source orbit space is compact and its target is Hausdorff, the induced bijection is a homeomorphism.

We must still pass from the covering manifold to its quotient by $\Gamma$. A representative of $f_{u,w}(xH)$ is $x'=r_uxw$, and
\begin{equation}\label{eq:deckonv}
 v'=(x')^{-1}u=w^{-1}x^{-1}r_uu=-w^{-1}v.
\end{equation}
On every nonsingleton block, this preserves the norm of $v_i$. On $V_0$, the identity $w|_{V_0}=-I_{V_0}$ gives $v'_0=v_0$. Thus $f$ fixes all the orbit invariants above. In particular, when $t=1$ it fixes the signed coordinate $z$. Consequently $f_{u,w}$ preserves each $C$-orbit, and
\[
 \Isom(\Gamma\backslash\Md)^0\backslash
       (\Gamma\backslash\Md)
 \cong C\backslash\Md.
\]
This proves the orbit-space assertions for the quotient itself.
\end{proof}

\textbf{Proof of Theorem~\ref{thm:main} (ii).}
If $t=0$, the cohomogeneity in Theorem~\ref{thm:orbits} is $q-1\ge2$. If $t>0$, the existence of at least one nonsingleton block gives $q-t\ge1$. Thus the full connected centralizer is never transitive. Since a finite covering and a finite effective kernel do not change orbit dimensions, Proposition~\ref{prop:cover} 
and Theorem \ref{thm:orbits} imply that $\Gamma\backslash\Md$ is not homogeneous.
\qed

We finish this section with the topological invariants of the quotient. Since $\Md$ is simply connected and $\Gamma\cong\Z_2$ is its deck group, we have
\begin{equation}\label{eq:fundamental}
 \pi_1(\Gamma\backslash\Md)\cong\Z_2,
 \qquad
 \chi(\Gamma\backslash\Md)=\frac12\chi(\Md).
\end{equation}

Let $a$ be the number of odd block dimensions, including singleton blocks, and put $b_i=\lfloor d_i/2\rfloor$ and $n=\lfloor N/2\rfloor$. Since $N=2\sum_i b_i+a$, we have
\begin{equation}\label{eq:ranks}
 \rank\SO(N)=\left\lfloor\frac N2\right\rfloor,
  \,\rank H=\sum_i b_i=\frac{N-a}{2},
  \,\rank\SO(N)-\rank H=\left\lfloor\frac a2\right\rfloor.
\end{equation}
The Hopf--Samelson theorem and the Weyl-group formula for the equal-rank case \cite{Borel,HopfSamelson} therefore apply directly.

\begin{proposition}\label{cor:euler}
Under the hypotheses of Theorem~\ref{thm:main},
\begin{equation}\label{eq:euler}
 \chi(\Gamma\backslash\Md)
 =\begin{cases}
 \displaystyle\frac{|W_{\SO(N)}|}{2\prod_i|W_{\SO(d_i)}|}
 =\displaystyle 2^{q-2}\frac{n!}{\prod_i b_i!},&a\le1,\\[6pt]
 0,&a\ge2.
 \end{cases}
\end{equation}
Here $W_K$ denotes the Weyl group of the compact connected Lie group $K$.
\end{proposition}
\begin{proof}
Note that $\chi(\Md)>0$ if and only if $\rank\SO(N)=\rank H$.
Equation~\eqref{eq:ranks} says that equality of ranks is equivalent to $a\le1$. The classical orders are
\[
 |W_{\SO(2b+1)}|=2^b b!,\qquad
 |W_{\SO(2b)}|=2^{b-1}b!\quad(b\ge1),
 \qquad |W_{\SO(1)}|=1.
\]
In particular, $\SO(2)$ has trivial Weyl group. Applying these formulas to the product isotropy and then using \eqref{eq:fundamental} gives the first expression in \eqref{eq:euler}.

If $a=0$, all $d_i=2b_i$ and $N=2n$, with $n=\sum_i b_i$. Thus
\[
\chi(\Md)= \frac{|W_{\SO(N)}|}{\prod_i|W_{\SO(d_i)}|}
 =\frac{2^{n-1}n!}{2^{\sum_i(b_i-1)}\prod_i b_i!}
 =2^{q-1}\frac{n!}{\prod_i b_i!}.
\]
If $a=1$, then $N=2n+1$ and precisely one block is odd, possibly a singleton. The odd block contributes $2^{b_i}b_i!$, including $b_i=0$, while each of the $q-1$ even blocks contributes $2^{b_i-1}b_i!$. Hence
\[
 \chi(\Md)=\frac{|W_{\SO(N)}|}{\prod_i|W_{\SO(d_i)}|}
 =\frac{2^n n!}{2^{n-(q-1)}\prod_i b_i!}
 =2^{q-1}\frac{n!}{\prod_i b_i!}.
\]
Dividing by two gives the factorial expression in both cases.
\end{proof}

\section{Stiefel and full flag counterexamples}\label{sec:examples}
The general construction becomes especially concrete for Stiefel manifolds and full flag manifolds. Clifford isometries of real Stiefel manifolds were already studied by C\'ampoli \cite{Campoli}.
\begin{corollary}\label{cor:stiefel}
Let $k\ge2$ be even, $d\ge2$, and $N=d+k\ge4$. The normal Stiefel manifold
\[
 V_k(\R^N)=\SO(N)/\SO(d)
\]
admits a free CW involution of displacement $\pi\sqrt{k(N-2)}$. Its quotient has fundamental group $\Z_2$, Euler characteristic zero, and connected isometry cohomogeneity one.
\end{corollary}
\begin{proof}
Take block dimensions $(d,1,\ldots,1)$, with $k$ singleton blocks. Then $q=k+1$ is odd and $t=k$. Theorem~\ref{thm:main} gives
\[
 \delta=\pi\sqrt{(N-2)(q-1)}=\pi\sqrt{k(N-2)},
 \qquad \cohom=q-t=1.
\]
At least two blocks are odd, so Corollary~\ref{cor:euler} gives zero Euler characteristic. The fundamental group is $\Z_2$ by \eqref{eq:fundamental}.
\end{proof}

\paragraph{Normal full flag manifolds.}
For a compact connected semisimple group $G$ and a maximal torus $T$, the normal full flag manifold is $G/T$ with the quotient metric induced by $-K_{\Lie(G)}$. If $p:\widetilde G\to G$ is a central covering and $\widetilde T=p^{-1}(T)$, then $\widetilde T$ is a maximal torus and
\[
 \widetilde G/\widetilde T\longrightarrow G/T,\qquad
 \widetilde x\widetilde T\longmapsto p(\widetilde x)T,
\]
is a diffeomorphism. Its differential preserves the Killing form, so it is an isometry for the corresponding normal metrics. Thus the following orthogonal models also represent full flags of simply connected compact groups.

\emph{Type $B_n$, with $n\ge2$ even.}
Take
\[
 \Md=\SO(2n+1)/T^n,
 \qquad \dd=(\underbrace{2,\ldots,2}_{n\text{ blocks}},1).
\]
Here $N=2n+1$, $q=n+1$ is odd, and $t=1$. Put
\[
 D=\diag(1,-1),\qquad
 b=\diag(\underbrace{D,\ldots,D}_{n\text{ blocks}},1),
 \qquad a_u=2uu^T-I_{2n+1}=-r_u.
\]
Because $n$ is even, $\det b=(-1)^n=1$. Each $D$ conjugates a plane rotation to its inverse, so $b\in N_{\SO(2n+1)}(T^n)$. The matrix $w=-b$ has determinant $-1$ on every two-plane block and on the singleton block. Since the two minus signs cancel,
\[
 f_{u,w}(xT^n)=a_uxbT^n,
 \qquad \delta_{f_{u,w}}=\pi\sqrt{n(2n-1)}.
\]
At least one block is even, so Proposition~\ref{prop:centralizer} gives
\[
 \Isom(\Gamma\backslash\Md)^0\cong\SO(2n),
 \qquad \cohom=n.
\]
There is exactly one odd block, and $|W_{B_n}|=2^n n!$, whereas the torus has trivial Weyl group. Consequently
\[
 \chi(\Md)=2^n n!,\qquad
 \chi(\Gamma\backslash\Md)=2^{n-1}n!.
\]

\emph{Type $D_n$, with $n\ge3$ odd.}
Now take
\[
 \Md=\SO(2n)/T^n,
 \qquad \dd=(\underbrace{2,\ldots,2}_{n\text{ blocks}}),
 \qquad w=\diag(D,\ldots,D).
\]
Here $q=n$ is odd, $t=0$ and $N=2n$. The matrix $w$ has determinant $-1$ and normalizes $T^n$ inside $\OO(2n)$. Thus
\[
f_{u,w}(xT^n)=r_uxwT^n,
 \qquad \delta_{f_{u,w}}=\pi\sqrt{(2n-2)(n-1)}
                         =\sqrt2\,\pi(n-1).
\]
The full connected quotient-isometry group and the orbit space are
\[
 \Isom(\Gamma\backslash\Md)^0\cong\SO(2n-1),
 \qquad \Delta^{n-1},
\]
so the cohomogeneity is $n-1$. Since $|W_{D_n}|=2^{n-1}n!$,
\[
 \chi(\Md)=2^{n-1}n!,\qquad
 \chi(\Gamma\backslash\Md)=2^{n-2}n!.
\]
The isomorphism $\Spin(6)\cong\SU(4)$ identifies $D_3=A_3$ and includes the normal full flag $\SU(4)/T^3$. Similarly, $B_2=C_2$ is the familiar identification $\Spin(5)\cong\Sp(2)$; it is not an additional construction. Table~\ref{tab:flags} summarizes these families.

\begin{table}[htbp]
\centering
\caption{Full-flag examples with positive Euler characteristic. Every quotient has fundamental group $\Z_2$.}
\label{tab:flags}
\small
\begin{tabular}{@{}llll@{}}
\toprule
Universal cover & Range & Quotient Euler characteristic & Cohomogeneity\\
\midrule
$\SO(2n+1)/T^n$ & $n\ge2$ even & $2^{n-1}n!$ & $n$\\
$\SO(2n)/T^n$ & $n\ge3$ odd & $2^{n-2}n!$ & $n-1$\\
\bottomrule
\end{tabular}
\end{table}

The product splitting of CW transformations on complete simply connected Riemannian manifolds belongs to the classical theory \cite[Theorem~3.1.2 and Corollary~3.1.3]{WolfSym}. For the next construction, the product distance and uniqueness of the de Rham decomposition suffice.

\begin{proposition}\label{prop:product}
Let $M_1,\ldots,M_r$ be compact, connected, simply connected, Riemannian irreducible manifolds, and let $\Gamma_i\le\Isom(M_i)$ be finite CW subgroups. Suppose that
\[
 \Gamma\le\prod_{i=1}^r\Gamma_i
\]
projects onto every $\Gamma_i$. With the product metric, $\Gamma$ is a free CW subgroup of $\Isom(\prod_iM_i)$, and
\[
 \Gamma\backslash\prod_iM_i\text{ is homogeneous}
 \quad\Longleftrightarrow\quad
 \Gamma_i\backslash M_i\text{ is homogeneous for every }i.
\]
\end{proposition}
\begin{proof}
For $\gamma=(\gamma_1,\ldots,\gamma_r)\in\Gamma$, the product distance satisfies
\[
 d\bigl((x_i),(\gamma_i(x_i))\bigr)^2
   =\sum_{i=1}^r d_i(x_i,\gamma_i(x_i))^2
   =\sum_{i=1}^r\delta_{\gamma_i}^2.
\]
Each summand is independent of $x_i$, so every $\gamma$ is CW. If $\gamma\ne1$, then some $\gamma_i\ne1$, and that component has strictly positive displacement. Thus $\gamma$ has no fixed point, proving freeness.

The factors are compact and simply connected, so there are no Euclidean factors in their de Rham decompositions. Connected isometries of their product preserve each irreducible factor; permutations of isometric factors belong to disconnected components. Consequently
\[
 \Isom\Bigl(\prod_iM_i\Bigr)^0=\prod_i\Isom(M_i)^0.
\]
For an element $(c_i)$ of this product, commuting with every $(\gamma_i)\in\Gamma$ is equivalent to
\[
 c_i\gamma_i=\gamma_ic_i
 \quad\text{for every }\gamma_i\text{ in the }i\text{-th projection of }\Gamma.
\]
The projection is $\Gamma_i$ by hypothesis. Taking identity components therefore yields
\[
 C_{\Isom(\prod_iM_i)}(\Gamma)^0
   =\prod_i C_{\Isom(M_i)}(\Gamma_i)^0.
\]
A product action is transitive precisely when each factor action is transitive. Applying Proposition~\ref{prop:cover} to the product and to each factor proves the equivalence.
\end{proof}
For example, if each $M_i$ is one of the irreducible oriented real flag manifolds of Theorem~\ref{thm:main}, one may use either the full product of the order-two deck groups or the diagonal order-two subgroup generated by $(f_1,\ldots,f_r)$. Both choices project onto every factor group and yield nonhomogeneous quotients. The diagonal involution has displacement
\[
 \delta_{(f_1,\ldots,f_r)}
 =\left(\sum_{i=1}^r\pi^2(N_i-2)(q_i-1)\right)^{1/2}.
\]
This provides further explicit examples without changing the normal metrics on the factors.

\section*{Declarations}

\noindent\textbf{Acknowledgments.}
Z. Yan is supported by the
Zhejiang Provincial Natural Science Foundation of China under Grant No. LMS25A010010. H. Zhang is
supported by NSFC (Grant No. 12301032), NSF of Jiangsu (Grant No. BK20230803) and Zhishan Scholars
Programs of Southeast University (Grant No. 2242024RCB0039).

\noindent\textbf{AI Disclosure.} Generative artificial intelligence tools, specifically ChatGPT and
WorkBuddy, were used during the preparation of this manuscript as auxiliary aids for language
editing, improving exposition and organization. All mathematical statements, proofs,
computations, references, and conclusions were independently verified by the authors, who take full
responsibility for the content of the manuscript.

\noindent\textbf{Data availability statement.} This manuscript has no associated data.

\noindent\textbf{Conflicts of interest.} The authors declare that they have no conflict of
interest.

\end{document}